# Response surface single loop reliability-based design optimization with higher-order reliability assessment

Rami Mansour[1] · Mårten Olsson[1]

**Abstract** Reliability-based design optimization (RBDO) aims at determination of the optimal design in the presence of uncertainty. The available Single-Loop approaches for RBDO are based on the First-Order Reliability Method (FORM) for the computation of the probability of failure, along with different approximations in order to avoid the expensive inner loop aiming at finding the Most Probable Point (MPP). However, the use of FORM in RBDO may not lead to sufficient accuracy depending on the degree of nonlinearity of the limit-state function. This is demonstrated for an extensively studied reliability-based design for vehicle crashworthiness problem solved in this paper, where all RBDO methods based on FORM strongly violates the probabilistic constraints. The Response Surface Single Loop (RSSL) method for RBDO is proposed based on the higher order probability computation for quadratic models previously presented by the authors. The RSSL-method bypasses the concept of an MPP and has high accuracy and efficiency. The method can solve problems with both constant and varying standard deviation of design variables and is particularly well suited for typical industrial applications where general quadratic response surface models can be used. If the quadratic response surface models of the deterministic constraints are valid in the whole region of interest, the method becomes a true single loop method with accuracy higher than traditional SORM. In other cases, quadratic response surface models are fitted to the deterministic constraints around the deterministic solution and the RBDO problem is solved using the proposed single loop method.

✉ Rami Mansour
  ramimans@kth.se

  Mårten Olsson
  mart@kth.se

[1] Department of Solid Mechanics, Royal Institute of Technology, Stockholm, Sweden

## 1 Introduction

In a deterministic design, statistical uncertainty in design variables and design parameters is not taken into account. This leads to a significantly simpler optimization problem which however can result in an unreliable design. If the variation of design variables and/or design parameters is taken into account, the feasible design space, limited by the probabilistic constraints, will be smaller since both feasibility and reliability needs to be fulfilled (Du and Chen 2004). The construction of the probabilistic constraints in Reliability-based Design Optimization (RBDO) (Valdebenito and Schuller 2010) has been subject to extensive research. Traditionally, a nonlinear optimization problem aiming at locating the Most Probable Point (MPP) followed by the First-Order Reliability Method (FORM) is used to compute the design reliability (Madsen et al. 1986). This leads to a nested double loop formulation which is computationally expensive. Therefore, research has been focusing on either improving the efficiency of the double loop, or decoupling the double loop method (Aoues and Chateauneuf 2010). In the first category the traditional Reliability Index Approach (RIA) (Tu et al. 1999) is improved by modifying the formulation of the probabilistic constraint leading to the Performance Measure Approach (PMA). In the second category, decoupled double loop approaches where the reliability analysis (inner loop) is separated from the optimization (outer loop) has been developed, see e.g. Du and Chen (2004). However the most efficient methods are single loop approaches for which the computational cost of RBDO is reduced to that of a deterministic optimization, but this necessitates accurate explicit expressions for the probabilistic constraints. It is however important to emphasize that for single-loop methods in RBDO, the accuracy in the optimal RBDO solution is generally reduced compared to double-loop methods or sequential decoupled methods since the MPP is not accurately located. However, the efficiency is dramatically increased. The aim of single-loop methods is therefore to achieve efficiency comparable to a deterministic problem but still maintaining good accuracy in engineering applications.

In this work, a Response Surface Single-Loop (RSSL) approach capable of solving RBDO problems with both constant and varying standard deviations is proposed. If the deterministic constraints are general quadratic functions, for instances if quadratic response surface models describing the system response in the whole region of interest are used, the method becomes a true single loop method. If, however, the deterministic constraints are not quadratic, the RBDO problem is solved using quadratic response surface models fitted around the deterministic solution.

This paper starts with an overview of the most common strategies in RBDO followed by the presentation of the RSSL method. Finally, benchmarks of typical RBDO problems are solved and the results compared to other methods with respect to accuracy and efficiency.

## 2 Reliability-based design optimization

### 2.1 Description of the problem

Reliability-based design optimization (RBDO) is an optimization problem aiming at locating the optimal design with the variations of the design variables and parameters being considered (Li et al. 2010). The RBDO problem can typically be formulated according to Shan and Wang (2008)

$$
\min_{\mathbf{d}, \boldsymbol{\mu}_\mathbf{x}} C(\mathbf{d}, \boldsymbol{\mu}_\mathbf{x}, \boldsymbol{\mu}_\mathbf{p})
$$
$$
\text{s.t Prob}[g_i(\mathbf{d}, \mathbf{x}, \mathbf{p}) < 0] \leq P_{fi,\text{all}}, \quad i = 1, ..., n
$$
$$
\mathbf{d}^{\text{Lower}} \leq \mathbf{d} \leq \mathbf{d}^{\text{Upper}}, \quad \boldsymbol{\mu}_\mathbf{x}^{\text{Lower}} \leq \boldsymbol{\mu}_\mathbf{x} \leq \boldsymbol{\mu}_\mathbf{x}^{\text{Upper}} \quad (1)
$$

where $\mathbf{d}$, $\mathbf{x}$, $\mathbf{p}$, $\boldsymbol{\mu}_\mathbf{x}$ and $\boldsymbol{\mu}_\mathbf{p}$ are the column vectors of deterministic design variables, random design variables, random parameters, mean value of $\mathbf{x}$ and mean value of $\mathbf{p}$, respectively. The objective function $C$ is to be minimized while the probability of violating the $i$th deterministic constraint $g_i$ should be smaller than the allowed probability of failure $P_{fi,\text{all}}$. The optimal design should be within the region given by the lower and upper bound on design variables, denoted with a superscript "Lower" and "Upper". The relations for the limits of the design region should be interpreted for each component of the vectors. The standard deviation of design variables $\boldsymbol{\sigma}_\mathbf{x}$ can either be constant or be allowed to vary.

### 2.2 Reliability assessment

The aim of the reliability assessment is to approximate the probability of failure defined in (1) or equivalently to estimate the generalized reliability index $\beta$ (Madsen et al. 1986) defined as

$$
\beta(\mathbf{d}, \mathbf{x}, \mathbf{p}) = -\Phi^{-1}(\text{Prob}[g(\mathbf{d}, \mathbf{x}, \mathbf{p}) < 0]). \quad (2)
$$

Index $i$ for the constraints has been dropped for simplicity. Observe that in the formulation of (2), a mapping of the probabilistic constraint has been performed.

The most widely used approximation of the generalized reliability index $\beta$ according to (2) is given by the First-Order Reliability Method (FORM) as

$$
\beta \approx \beta^{\text{HL}} = \left\| \mathbf{v}_N^{\text{MPP}} \right\|, \quad (3)
$$

where an isoprobabilistic transformation from the original design space $\mathbf{v} = [\mathbf{x}^T, \mathbf{p}^T]^T$ to a standardized normal variable space $\mathbf{v}_N = [\mathbf{x}_N{}^T, \mathbf{p}_N{}^T]^T$ has been performed. In the Reliability-Index Approach (RIA), the point $\mathbf{v}_N^{\text{MPP}}$ is called the Most Probable Point (MPP) and is the shortest distance to the failure surface $g_N = 0$ in the $\mathbf{v}_N$-space, and is therefore found by solving the so called MPP problem (Hasofer and Lind 1974)

$$
\begin{cases} \beta^{\text{HL}} = \min_{\mathbf{v}_N} \|\mathbf{v}_N\| \\ \text{s.t. } g_N = 0 \end{cases} \quad (4)
$$

Another method is the Performance Measure Approach (PMA) (Tu et al. 1999) in which the probabilistic constraints in (1) are approximated by

$$
g(\mathbf{d}, \mathbf{v}^{\text{MPP}}) = 0, \quad (5)
$$

where $\mathbf{v}^{\text{MPP}}$ is the point in the original design space corresponding to the inverse MPP, $\mathbf{v}_N^{\text{MPP}}$, in the normalized variable space. As opposed to the Reliability Index Approach, the point $\mathbf{v}_N^{\text{MPP}}$ is found by solving the inverse MPP problem defined as

$$
\begin{cases} \min_{\mathbf{v}_N} g_N(\mathbf{d}, \mathbf{v}_N) \\ \text{s.t } \|\mathbf{v}_N\| = \beta_d \end{cases} \quad (6)
$$

where $\beta_d = \Phi^{-1}(r_d)$ is the desired reliability index and $r_d = 1 - P_f$ is the desired reliability. It should be noted that FORM is based on a linearisation of the failure surface at the MPP and its accuracy decreases with increased non-linearity of the failure surface in the $\mathbf{v}_N$-space. A higher order asymptotic expansion for the probability content Prob $[g(\mathbf{d}, \mathbf{x}, \mathbf{p}) < 0]$ has been derived by Breitung (Breitung 1984), thus yielding a higher order approximation of the generalized reliability index according to

$$\beta \approx -\Phi^{-1}\left[\Phi\left(-\beta^{HL}\right)\prod_i\left(1-\beta^{HL}\rho_i^{MPP}\right)^{-1/2}\right], \quad \beta^{HL} \to \infty \quad (7)$$

where $\rho_i^{MPP}$ are the main curvatures of the failure surface at the MPP. This Second-Order Reliability Method (SORM) is based on a hyper-parabolic approximation of the failure surface at the MPP. However, for small to moderate $\beta^{HL}$, the Breitungs' formula yields less accurate results (Madsen et al. 1986) and therefore other SORM approximations have been presented, see e.g. by Tvedt (1984). It has however been shown in Mansour and Olsson (2014) and Zhao and Ono (1999) that the most commonly used SORMs have one or more of the following drawbacks: They do not work well for small curvatures at the MPP and/or large number of random variables, neither do they work well for negative or strongly uneven curvatures at the MPP. The SORM expressions may even contain singularities.

### 2.3 Double loop method

The double loop method in RBDO is the most direct approach (Valdebenito and Schuller 2010) but also numerically expensive, since the structural reliability is estimated for each set of design variables evaluated by the optimization algorithm. In FORM and SORM, the reliability assessment is in itself an optimization problem. Therefore, the RBDO problem consists of two nested loops; an outer loop for the minimization of the objective function and an inner loop for evaluating the constraint, i.e. for the reliability assessment. Typically, the number of function evaluations using a double loop optimization is large, and different RBDO methods aims at easing the computation cost especially for problems with a large number of constraints (Dersjo and Olsson 2011). For a crashworthiness problem (Gu et al. 2001) widely used in comparing different RBDO methods, see e.g. Shan and Wang (2008) and Li et al. (2010), consisting of 10 constraints, 9 random design variables of which 2 are discrete, and 2 random parameters, the total number of function evaluations (Du and Chen 2004) for achieving a reliability level of $r_d = 0.9$ was 3,324 using the Performance Measure Approach (PMA) and 26,984 using the Reliability Index Approach (RIA). Most of the function evaluations are used for the reliability assessment. The flowchart of a typical double loop RBDO procedure is provided in Fig. 1.

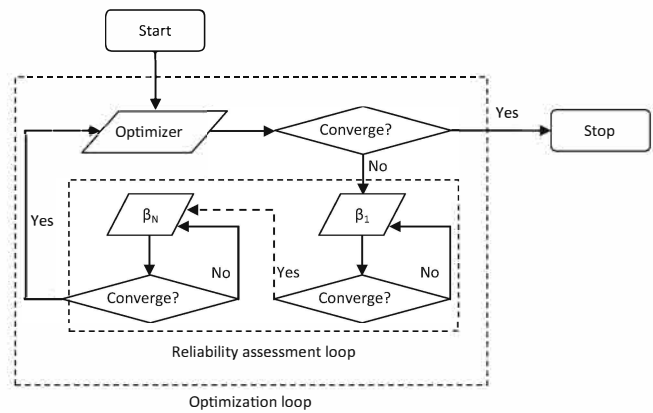

Fig. 1 Flowchart of a double-loop RBDO method

### 2.4 Decoupling method

In order to avoid the expensive double loop, a decoupled approach aiming at separating the reliability assessment from the outer optimization loop is used. One method is to perform a linear sequential optimization where a linear approximation of the reliability index is constructed using sensitivity information, i.e.

$$\beta(\boldsymbol{\mu}_\mathbf{x}) = \beta\left(\boldsymbol{\mu}_\mathbf{x}^{(k)}\right) + \nabla_{\boldsymbol{\mu}_\mathbf{x}}\beta\big|_{\boldsymbol{\mu}_\mathbf{x}=\boldsymbol{\mu}_\mathbf{x}^{(k)}}\left(\boldsymbol{\mu}_\mathbf{x} - \boldsymbol{\mu}_\mathbf{x}^{(k)}\right) \quad (8)$$

where $\boldsymbol{\mu}_\mathbf{x}^{(k)}$ is the $k$-th candidate optimal solution. A new candidate solution $\boldsymbol{\mu}_\mathbf{x}^{(k+1)}$ is found based on the linear approximation of the reliability index according to (8) which uses the sensitivities from the previous solution $\boldsymbol{\mu}_\mathbf{x}^{(k)}$. The process is repeated until convergence. The flowchart of this method is presented in Fig. 2. As can be seen, the major difference between the double loop and the decoupled loop is that the reliabilities do not need to be computed for each set of design variables evaluated by the optimization algorithm; instead they are only updated after the convergence of the linear optimization problem. Another decoupled loop approach is the Sequential Optimization and Reliability Assessment (SORA) (Du and Chen 2004) method, which approximate the probabilistic constraints by shifting the deterministic constraints towards the feasible region. An approximate optimal design is then found and the process is repeated until convergence.

### 2.5 Single loop method

In a single loop RBDO, the inner reliability assessment in the double loop method is avoided. Ideally, each probabilistic constraint is approximated by a function $g^*$ of the design variables thus reducing RBDO to a deterministic optimization problem, see the flowchart according to Fig. 3. The concept of reliable design space (RDS) (Shan and Wang 2008) is a complete single loop method as described in the flowchart. However RDS is based on approximation of FORM. This may lead to inaccuracies depending on the non-linearity of the problem. To the best of the authors' knowledge, all single-loop methods presented in literature are based on FORM. Furthermore, most single loop methods are not complete single-loops as described in the flowchart Fig. 3. They either increases the dimension of the problem, see e.g. Agarwal et al. (2007), or involves iterative processes in the evaluation of derivatives and/or constraints, see e.g. the Single-Loop Single Variable (SLSV) method (Chen et al. 1997), the improvement of SLSV presented by Liang et al. (2008) and the Traditional Approximation Method (TAM) (Yang and Gu 2004).

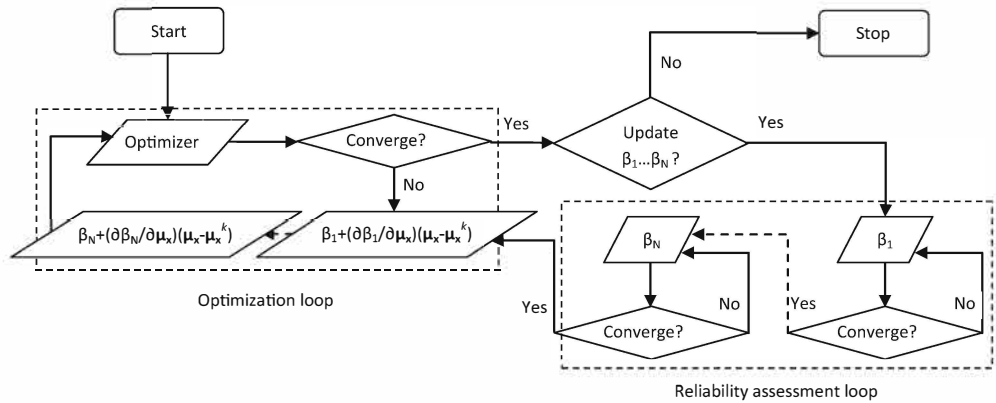

**Fig. 2** Flowchart of a decoupled-loop RBDO method

### 2.6 The use of response surface models

The use of response surface models in RBDO is necessary when analytical model is not available to describe the response of the system. This is the case in most industrial applications in which the Finite Element Method (FEM) is used. A review of four widely used response surface models is given in Jin et al. (2001): second-order Polynomial Regression (PR), the Kriging Method, Multivariate Adaptive Regression Splines and Radial Basis Functions.

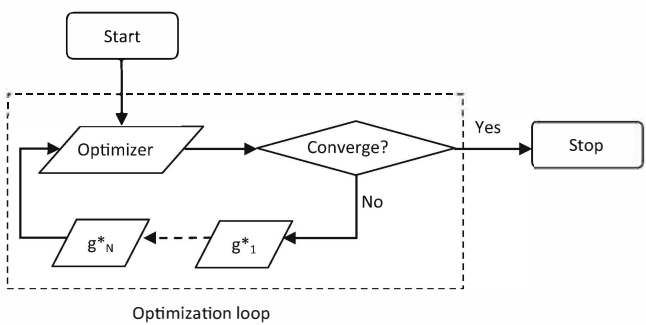

**Fig. 3** Flowchart of a single-loop RBDO method

They were compared with respect to Accuracy and Robustness, Efficiency, Transparency and Simplicity; where Transparency was defined as the ability to obtain the contribution of each input factor and the interaction between them. It was proposed that when constructing response surface models, PR should be implemented first to see if a reasonable fit can be obtained. This conclusion was based on the fact that second order polynomial regressions are easy to use, very accurate for low order non-linearity and have good accuracy. Furthermore the robustness of PR in problems with noise was superior, a conclusion which was previously made by e.g. Guinta et al. (1994), which noted that the smoothing capability of polynomial regressions allow quick convergence of noisy functions in optimization problems.

The practicability of using second-order polynomials in simulation based response surface design is enhanced by the implementation of effective built in functions in commercial FE codes. In ANSYS (2009) a central composite design (CCD) sampling specifically adapted to second-order models aiming at reducing the number of runs at a reasonable level while maintaining good accuracy is used. The CCD scheme implemented consists of a central point, $2N$ axis point and $2^{N-f}$ factorial points located at the corner of an $N$-dimensional hypercube. Here, $N$ is the number of random input variables and $f$ is the fraction of factorial part of the central composite design. The value of $f$ is determined automatically by ANSYS and increases gradually with increased number of random input variables, see Table 1. The CCD together with the Box-Behnken Design (BBD) are the most popular designs of experiments (DOE) for fitting a second-order polynomial response surface model (Montgomery 2001). The BBD is also implemented in ANSYS and consist of a central point and the midpoints of each edge of an $N$-dimensional hypercube and is usually very efficient in terms of the number of required runs, see Table 1. It should however be noted that BBD was originally proposed for 3-7, 9-12 or 16 factors. For other number of factors, the design is constructed in similar ways. It should be noted that the superiority of polynomial regressions in many levels, see Jin et al. (2001) and Guinta et al. (1994), is mostly related to quadratic response surface models, since higher order polynomials may cause instabilities (Barton 1992).

**Table 1** Number of sample points (simulation loops) required for a central composite design (CCD) and a Box-Behnken design (BBD) as function of the number of input variables and number of coefficients for a general quadratic function in ANSYS (2009)

| Number of variables | Number of coefficient | Central composite | | Box-Behnken |
|---|---|---|---|---|
| | | Factorial number | Sample points | Sample points |
| 1 | 3 | N/A | N/A | N/A |
| 2 | 6 | 0 | 9 | N/A |
| 3 | 10 | 0 | 15 | 13 |
| 4 | 15 | 0 | 25 | 25 |
| 5 | 21 | 1 | 27 | 41 |
| 6 | 28 | 1 | 45 | 49 |
| 7 | 36 | 1 | 79 | 57 |
| 8 | 45 | 2 | 81 | 65 |
| 9 | 55 | 2 | 147 | 121 |
| 10 | 66 | 3 | 149 | 161 |
| 11 | 78 | 4 | 151 | 177 |
| 12 | 91 | 4 | 281 | 193 |

## 3 Response Surface Single Loop (RSSL) method

In this work, the Response Surface Single Loop (RSSL) method capable of handling both constant and varying standard deviation of design variables is proposed. A flowchart of the method is presented in Fig. 4. In RSSL, the deterministic constraints are approximated by quadratic response surface models around the deterministic solution. Thereafter, the RBDO problem is solved using a single-loop approach where the reliability assessment is based on closed-form second-order probability formulas which takes into account all components of the Hessian of the quadratic model (Mansour and Olsson 2014).

### 3.1 Essence of the method

In RSSL, the deterministic constraints $g_i(\mathbf{d}, \mathbf{x}, \mathbf{p})$, recalling the RBDO problem according to (1), are approximated by general quadratic response surface models $Q_i$ in $\mathbf{d}$, $\mathbf{x}$ and $\mathbf{p}$ around the deterministic solution. Introducing

$$\mathbf{z} = \begin{bmatrix} \mathbf{d}^T & \mathbf{x}^T & \mathbf{p}^T \end{bmatrix}^T, \quad (9)$$

the response surface model can be written as

$$Q_i(\mathbf{z}) = \mathbf{z}^T \mathbf{A} \mathbf{z} + \mathbf{k}^T \mathbf{z} + c. \quad (10)$$

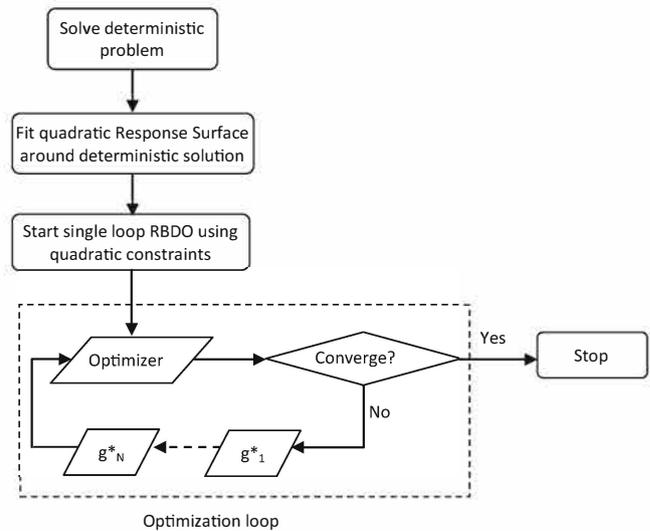

**Fig. 4** Flowchart of the proposed Response Surface Single-Loop (RSSL) method

If needed $\mathbf{A}$, $\mathbf{k}$ and $c$ may be replaced by a function of the deterministic design variables $\mathbf{d}$ throughout this paper. The quadratic dependence in $\mathbf{x}$ and $\mathbf{p}$ must, however, be retained. The deterministic design variables $\mathbf{d}$ are considered as special random variables with zero standard deviation and mean values $d_i$.

The vector $\mathbf{z}$ may consist of correlated non-normal variables with probability density function (PDF) $f_{zi}(\mu_{zi}, \sigma_{zi})$, cumulative density function (CDF) $F_{zi}(\mu_{zi}, \sigma_{zi})$ and matrix of correlation coefficients $\rho_{z_i,z_j}$ according to

$$\mathbf{C} = \begin{bmatrix} 1 & \rho_{z_1,z_2} & \cdots & \rho_{z_1,z_n} \\ \rho_{z_2,z_1} & 1 & \cdots & \rho_{z_2,z_n} \\ \vdots & \vdots & \vdots & \vdots \\ \rho_{z_n,z_1} & \rho_{z_n,z_2} & \cdots & 1 \end{bmatrix}. \quad (11)$$

The vector $\mathbf{z}$ is expressed in terms of standard normal and uncorrelated variables $\mathbf{z_N}$ according to (Haldar and Mahadevan 2000)

$$\mathbf{z} = \mathbf{S}_{\mathbf{z},eq} \mathbf{T} \mathbf{D} \mathbf{z}_N + \boldsymbol{\mu}_{\mathbf{z},eq}, \quad (12)$$

where $\mathbf{T}$ is an orthogonal transformation matrix with columns consisting of the normalized eigenvectors of $\mathbf{C}$, $\mathbf{D}$ is a diagonal matrix consisting of the square root of the corresponding eigenvalues $\lambda_i$ of $\mathbf{C}$

$$\mathbf{D} = \mathrm{diag}\left[\sqrt{\lambda_1}, \sqrt{\lambda_2}, ..., \sqrt{\lambda_n}\right], \quad (13)$$

$\mathbf{S}_{\mathbf{z},eq}$ is a diagonal matrix with the equivalent normal standard deviations $\sigma_{zi,eq}$

$$\mathbf{S}_{\mathbf{z},eq} = \mathrm{diag}\left[\sigma_{z1,eq}, \sigma_{z2,eq}, ..., \sigma_{zn,eq}\right] \quad (14)$$

and $\boldsymbol{\mu}_{\mathbf{z},eq}$ is the vector of equivalent normal means $\mu_{zi,eq}$. The equivalent normal mean values $\mu_{zi,eq}$ and equivalent normal standard deviations $\sigma_{zi,eq}$ are computed using the formulas for equivalent normalization (Choi et al. 2007) according to

$$\begin{cases} \sigma_{zi,\text{eq}} = \frac{\phi(\Phi^{-1}[F_{zi}(\mu_{zi},\sigma_{zi})])}{f_{zi}(\mu_{zi},\sigma_{zi})} \\ \mu_{zi,\text{eq}} = \mu_{zi} - \Phi^{-1}\left[F_{zi}(\mu_{zi},\sigma_{zi})\right]\sigma_{zi,\text{eq}} \end{cases}, \quad (15)$$

where $\Phi$ and $\phi$ are the standard normal cumulative density function (CDF) and probability density function (PDF), respectively. For normal random variables $\mathbf{z}$, $\mu_{zi,\text{eq}}$ and $\sigma_{zi,\text{eq}}$ are simply given by

$$\begin{cases} \sigma_{zi,\text{eq}} = \sigma_{zi} \\ \mu_{zi,\text{eq}} = \mu_{zi} \end{cases}, \quad (16)$$

and for the deterministic components of $\mathbf{z}$, see (9), (15) is given by

$$\begin{cases} \sigma_{di,\text{eq}} = 0 \\ \mu_{di,\text{eq}} = d_i \end{cases}. \quad (17)$$

For uncorrelated variables, $\mathbf{T} = \mathbf{D} = \mathbf{I}$, where $\mathbf{I}$ is the identity matrix.

The response surface can thus be transformed using (12) according to

$$Q_{Ni}(\mathbf{z}_N) = Q_i(\mathbf{z}) \quad (18)$$

where

$$Q_{Ni}(\mathbf{z}_N) = \mathbf{z}_N^T \mathbf{A}' \mathbf{z}_N + \mathbf{k}'^T \mathbf{z}_N + c', \quad (19)$$

and

$$\begin{cases} \mathbf{A}' = \mathbf{D}\mathbf{T}^T \mathbf{S}_{z,\text{eq}} \mathbf{A} \mathbf{S}_{z,\text{eq}} \mathbf{T} \mathbf{D} \\ \mathbf{k}'^T = \mathbf{k}^T \mathbf{S}_{z,\text{eq}} \mathbf{T} \mathbf{D} + 2\boldsymbol{\mu}_{z,\text{eq}}^T \mathbf{A} \mathbf{S}_{z,\text{eq}} \mathbf{T} \mathbf{D} \\ c' = c + \boldsymbol{\mu}_{z,\text{eq}}^T \mathbf{A} \boldsymbol{\mu}_{z,\text{eq}} + \mathbf{k}^T \boldsymbol{\mu}_{z,\text{eq}} \end{cases}. \quad (20)$$

Observe that the components of $\mathbf{S}_{z,\text{eq}}$ corresponding to the deterministic variables is simply equal to zero and the corresponding components of $\boldsymbol{\mu}_{z,\text{eq}}$ is equal to $d_i$, according to (17).

For RBDO problems with varying standard deviation of design variables, typically the standard deviation is given by $\sigma_{xi} = t_i \mu_{xi}$, where $t_i$ is a constant, see e.g. (Yin and Chen 2006); which after insertion in the expression for $\mathbf{S}_{z,\text{eq}}$ and $\boldsymbol{\mu}_{z,\text{eq}}$ eliminates the dependence on $\sigma_{xi}$ in (20). Therefore, RSSL can also handle varying standard deviation of design variables with the same formulation. The probability of violating the deterministic constraint in (1) is now approximated using the quadratic response surface according to

$$P_f = \text{Prob}\left[\mathbf{z}_N^T \mathbf{A}' \mathbf{z}_N + \mathbf{k}'^T \mathbf{z}_N + c' < 0\right], \quad (21)$$

for which Mansour and Olsson (2014) proposed two closed form expressions depending on the sign of the eigenvalues of $\mathbf{A}'$. These formulas are restated here. Denote by $\gamma_j$ the eigenvalues of $\mathbf{A}'$ and by $\mathbf{P}$ the orthogonal transformation matrix with columns consisting of the corresponding normalized eigenvectors of $\mathbf{A}'$. Define

$$\begin{cases} \bar{\mathbf{k}}^T = \mathbf{k}'^T \mathbf{P} \\ m_r = \sum_j \left(\gamma_j^r + \frac{r}{4}\gamma_j^{r-2}\bar{k}_j^2\right), r = 1, 2, 3, 4 \end{cases}. \quad (22)$$

For the case where $\gamma_j$ have different signs, the probability of failure is approximated by the closed-form expression

$$P_f = \Phi(\kappa_1) - \varphi(\kappa_1)\left[\begin{array}{c} \frac{\sqrt{2}H_2(\kappa_1)}{3} \frac{m_3}{m_2\sqrt{m_2}} \\ + \frac{H_5(\kappa_1)}{9} \frac{m_3^2}{m_2^3} + \frac{H_3(\kappa_1)}{2} \frac{m_4}{m_2^2} \end{array}\right], \quad (23)$$

where $H_i$ is the $i$:th probabilists' Hermite polynomials and

$$\kappa_1 = -\frac{c' + \sum_j \gamma_j}{\sqrt{\sum_j \left(2\gamma_j^2 + \bar{k}_j^2\right)}}. \quad (24)$$

For a linear approximation (FORM), the probability of failure is reduced to the exact value $P_f = \Phi(-\beta)$, since then $\gamma_j = 0$ for all $j$, $\bar{\mathbf{k}} = \begin{bmatrix} 0 & \cdots & 0 & -1 \end{bmatrix}^T$ and $\kappa_1 = -\beta$.

For the case where all $\gamma_j$ have the same sign, which corresponds to elliptic limit-states, the probability of failure is approximated by

$$P_f = \begin{cases} P, & \text{sign}(\gamma_j) \cdot h > 0 \\ 1 - P, & \text{sign}(\gamma_j) \cdot h < 0 \end{cases}, \quad (25)$$

where

$$P = \Phi(\kappa_2) - \varphi(\kappa_2)\left[\begin{array}{c} H_3(\kappa_2)\left(\frac{m_4^2}{2m_2^2} - \frac{20m_3^2}{27m_2^3}\right) \\ + \frac{2m_3}{9m_1 m_2} \\ + H_1(\kappa_2)\left(-\frac{2m_3^2}{3m_2^3} + \frac{2m_3}{3m_1 m_2}\right) \end{array}\right] \quad (26)$$

and

$$\begin{cases} \kappa_2 = \frac{|m_1|}{\sqrt{2h^2 m_2}}\left[\left(\left|\frac{q_0}{m_1}\right|\right)^h - 1 - \frac{h(h-1)m_2}{m_1^2}\right] \\ h = 1 - \frac{2m_1 m_3}{3m_2^2} \\ q_0 = \sum_j \frac{\bar{k}_j^2}{4\gamma_j} - c' \end{cases}. \quad (27)$$

In order to avoid singularities in (27), if any $\gamma_j = 0$, it should be replaced by $\gamma_j = \varepsilon$ when all $\gamma_j \geq 0$ and by $\gamma_j = -\varepsilon$ when all $\gamma_j \leq 0$. Here $\varepsilon$ is a sufficiently small number, in this work $\varepsilon = 10^{-7}$.

Since the expressions for the probability of failure can be applied directly given (20) without any MPP-search, they are used to form the probabilistic constraints and the RBDO problem according to (1) is solved in a single loop with standard optimization tools. Therefore, since the probabilistic constraints are explicitly written as functions of the design variables which typically are the mean values of the random variables, finite difference methods can be used directly on the constraints to compute sensitivities. Furthermore, the probabilistic constraints are analytically given based on the fitting parameters of the response surface and therefore no additional function evaluation of the deterministic constraints is performed when computing sensitivities. This is an important advantage compared to MPP methods, where the probabilistic constraints vary with the MPP. Therefore, computing the sensitivities of the probability of failure with respect to mean values using a finite difference method is expensive for MPP-based RBDO methods. It necessitates the expensive computation of the MPP at a finite distance from the actual point. For these methods, analytical sensitivities are necessary and crucial for efficiency in RBDO. The probabilistic optimization problem in RSSL is therefore reduced to a typical deterministic problem since constraints are analytically given and therefore standard optimization tools such as the widely used MatLabs "fmincon" built-in function, see e.g. Shan and Wang (2008), can be applied.

The fitting of the quadratic response surface model around the deterministic solution according to (10), is made in the original $\mathbf{z}$-space using a Box-Behnken design, see Table 1. Although the box domain is made in the original design space with arbitrarily distributed and correlated design variables, the size of the box domain need to be decided in the standardized uncorrelated space $\mathbf{z}_N$, and then transformed back to the original space. The size of the box domain depends on the desired reliability level and a constant $C_R$, according to, see Lee et al. (2011),

$$\mathbf{z}_N : \pm C_R \beta_d. \quad (28)$$

The value of $C_R$ should be chosen larger than unity for the design variables so that the constraints in $\mathbf{z}_N$-space evaluated at the optimal mean value lies within the circle with radius $C_R$ and therefore including the MPP at the predicted optimal solution, see Fig. 5 and Lee et al. (2011). A relatively large value of $C_R = 1.4$ is used here for the design variables, to assure that the response surface is valid in a larger domain since no update of the response surface is allowed in order to qualify the method as a single-loop method. The used intervals around the deterministic solution are

$$\begin{cases} \mathbf{x} : & \pm 1.4 \beta_d \sigma_{x_i,\text{eq}} \\ \mathbf{p} : & \pm \beta_d \sigma_{p_i,\text{eq}} \\ \mathbf{d} : & \pm 1.4 \beta_d \frac{d_i}{10} \end{cases}, \quad (29)$$

where the equivalent normal standard deviation $\sigma_{zi,\text{eq}}$ is computed at the deterministic solution $\mu_{zi,\text{det}}$ according to

$$\sigma_{zi,\text{eq}} = \frac{\phi\left(\Phi^{-1}\left[F_{zi}\left(\mu_{zi,\text{det}}, \sigma_{zi}\right)\right]\right)}{f_{zi}\left(\mu_{zi,\text{det}}, \sigma_{zi}\right)}, \quad (30)$$

where $C_R = 1$ for $\mathbf{p}_N$ and $C_R = 1.4$ for $\mathbf{x}_N$ in (28). It should also be noted that the size of the box domain for the deterministic variables $\mathbf{d}$ can be chosen more freely. Furthermore, since $\mu_p$ does not change during the search for the probabilistic optimum, the $x_N$-component of the MPP may be expected to be further away from the origin compared to $p_N$-component, see Fig. 5b). Therefore a smaller value of $C_R = 1$ for the design parameters is used compared to the larger value of $C_R = 1.4$ used for the design variables.

### 3.2 Illustration of the probabilistic constraints

#### 3.2.1 Uncorrelated normally distributed random variables with constant standard deviation

For illustration purposes the construction of the probabilistic constraint given the following deterministic constraint

$$g(x_1, p_1) = \frac{(x_1 + p_1 - 5)^2}{30} + \frac{(x_1 - p_1 - 12)^2}{120} - 1 \quad (31)$$

is demonstrated. The random design variable is normal with mean value $0 \leq \mu_{x1} \leq 15$ and standard deviation $\sigma_{x1} = 0.3$. The random design parameter $p_1$ is also normal with mean value $\mu_{p1} = 3.4$ and $\sigma_{p1} = 0.3$. Both random

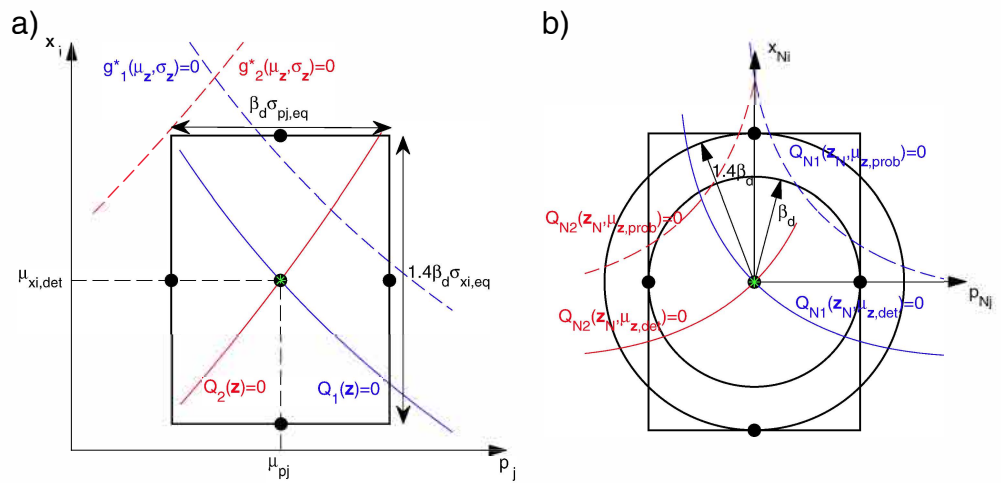

**Fig. 5** A Box-Behnken Design of Experiments in RSSL in **a** original space and **b** standard normalized space. The fitted quadratic functions in original space are denoted by $Q_i(\mathbf{z})$ and the RSSL-probabilistic constraints by $g_i^*(\mu_\mathbf{z}, \sigma_\mathbf{z})$. The transformed quadratic functions evaluated at the deterministic solution and at the probabilistic solution are denoted by $Q_{Ni}(\mathbf{z}_N, \mu_{\mathbf{z},\text{det}})$ and $Q_{Ni}(\mathbf{z}_N, \mu_{\mathbf{z},\text{prob}})$, respectively

variables are assumed uncorrelated. Since (31) is already quadratic, it can be written according to (10) as

$$g(\mathbf{z}) = Q(\mathbf{z}) = \mathbf{z}^T \mathbf{A} \mathbf{z} + \mathbf{k}^T \mathbf{z} + c, \quad (32)$$

where

$$\mathbf{z} = \begin{bmatrix} z_1 & z_2 \end{bmatrix}^T = \begin{bmatrix} x_1 & p_1 \end{bmatrix}^T \quad (33)$$

and

$$\begin{cases} \mathbf{A} = \begin{bmatrix} 1/24 & 1/40 \\ 1/40 & 1/24 \end{bmatrix} \\ \mathbf{k} = \begin{bmatrix} -8/15 & -2/15 \end{bmatrix}^T \\ c = 31/30 \end{cases} \quad (34)$$

Since the random variables $z_1$ and $z_2$ are assumed uncorrelated, the matrix of correlation coefficients according to (11) is reduced $\mathbf{C} = \mathbf{I}$. The orthogonal transformation matrix $\mathbf{T}$ and the diagonal matrix $\mathbf{D}$ used in (20), see (12) and (13), are therefore given by

$$\mathbf{T} = \mathbf{D} = \begin{bmatrix} 1 & 0 \\ 0 & 1 \end{bmatrix}. \quad (35)$$

Since the variables are normally distributed, the vector of equivalent means $\boldsymbol{\mu}_{\mathbf{z},\text{eq}}$ and the matrix of equivalent standard deviations $\mathbf{S}_{\mathbf{z},\text{eq}}$ according to (14), are given by

$$\boldsymbol{\mu}_{\mathbf{z},\text{eq}} = \begin{bmatrix} \mu_{x1} & \mu_{p1} \end{bmatrix}^T = \begin{bmatrix} \mu_{x1} & 3.4 \end{bmatrix}^T \quad (36)$$

and

$$\mathbf{S}_{\mathbf{z},\text{eq}} = \begin{bmatrix} \sigma_{x1} & 0 \\ 0 & \sigma_{p1} \end{bmatrix} = \begin{bmatrix} 0.3 & 0 \\ 0 & 0.3 \end{bmatrix}, \quad (37)$$

respectively. Therefore $\mathbf{A}'$, $\mathbf{k}'$ and $c'$ in (20) are given by:

$$\begin{cases} \mathbf{A}' = \begin{bmatrix} 0.0037 & 0.0022 \\ 0.0022 & 0.0037 \end{bmatrix} \\ \mathbf{k}' = \begin{bmatrix} \mu_{x1}/40 - 0.109 \\ 3\mu_{x1}/200 + 0.045 \end{bmatrix} \\ c' = \mu_{x1}^2/24 - 0.3633\mu_{x1} + 1.0617 \end{cases} \quad (38)$$

The orthogonal transformation matrix with columns consisting of the normalized eigenvectors of $\mathbf{A}'$ is given by

$$\mathbf{P} = \frac{1}{\sqrt{2}} \begin{bmatrix} -1 & 1 \\ 1 & 1 \end{bmatrix} \quad (39)$$

and the corresponding eigenvalues of $\mathbf{A}'$ by

$$\begin{cases} \gamma_1 = 0.0015 \\ \gamma_2 = 0.0060 \end{cases}. \quad (40)$$

Therefore, following (22),

$$\bar{\mathbf{k}} = \sqrt{2} \begin{bmatrix} 0.077 - 0.005\mu_{x1} \\ 0.02\mu_{x1} - 0.032 \end{bmatrix} \quad (41)$$

and

$$\begin{cases} m_1 = 0.0417\mu_{x1}^2 - 0.363\mu_{x1} + 2.07 \\ m_2 = 4.25 \cdot 10^{-4}\mu_{x1}^2 - 0.00205\mu_{x1} + 0.00699 \\ m_3 = 3.66 \cdot 10^{-6}\mu_{x1}^2 - 1.33 \cdot 10^{-5}\mu_{x1} + 2.28 \cdot 10^{-5} \\ m_4 = 2.89 \cdot 10^{-8}\mu_{x1}^2 - 9.56 \cdot 10^{-8}\mu_{x1} + 1.02 \cdot 10^{-7} \end{cases} \quad (42)$$

It should be noted that $\bar{\mathbf{k}}$ and $m_r$ are functions of the mean value of the design variables. Since the eigenvalues are of the same sign (elliptic $g(\mathbf{z})$), the probability of failure is given by (25). It can be easily shown that for $0 \leq \mu_{x1} \leq 15$

$$h = 1 - 2m_1 m_2 / 3m_2^2 > 0, \quad (43)$$

see (27). Therefore $\text{sign}(\gamma_j) \cdot h > 0$, and the probability of failure is given by $P_f = P$ where $P$ is given by (26). Computing $q_0$ in (27) using (38), (40) and (41) yields

$$q_0 = \sum_j \frac{\bar{k}_j^2}{4\gamma_j} - c' = 1. \quad (44)$$

Using that $m_1 > 0$ for $0 \leq \mu_{x1} \leq 15$ and inserting (43) and (44) into the expression of $\kappa_2$ in (27) yields

$$\kappa_2 = \frac{m_1^{(2m_1 m_3/3m_2^2)} - \frac{4m_1 m_3^2}{9m_2^3} + \frac{2m_3}{3m_2} - m_1}{\sqrt{2m_2 \left(\frac{2m_1 m_3}{3m_2^2} - 1\right)^2}}. \quad (45)$$

Assuming that the desired reliability is $\beta_d = 3$, the allowed probability of failure is given by $P_{f,\text{all}} = \Phi(-\beta_d) = 0.0013$ and the probabilistic constraint by $g^*(\mu_{x1}) = P_{f,\text{all}} - P_f$:

$$g^*(\mu_{x1}) = P_{f,\text{all}} - \Phi(\kappa_2) \\ + \varphi(\kappa_2) \begin{bmatrix} H_3(\kappa_2) \left(\frac{m_4^2}{2m_2^2} - \frac{20m_3^2}{27m_2^3} + \frac{2m_3}{9m_1 m_2}\right) \\ + H_1(\kappa_2) \left(-\frac{2m_3^2}{3m_2^3} + \frac{2m_3}{3m_1 m_2}\right) \end{bmatrix}. \quad (46)$$

In (45), the expression for $h$ in (43) and $q_0 = 1$ in (44) has been used. It should be noted that $g^*(\mu_{x1})$ is a function of $\mu_{x1}$ only since $m_1$, $m_2$, $m_3$ and $m_4$ in (42) are functions of $\mu_{x1}$. Numerical sensitivities are therefore easily computed using (46) in RBDO applications without any function evaluation of the deterministic constraint in (31).

The deterministic constraint $g(\mathbf{z})$ as well as the probability of failure as a function of the mean value of the design variables $\mu_{x1}$ for $\mu_{z2} = \mu_{p1} = 3.4$ is shown in Fig. 6. As can be seen the feasible region for the mean value of the design variable is $0 \leq \mu_{x1} \leq 3.86$ and $5.93 \leq \mu_{x1} \leq 15$. In the non-feasible region $3.86 < \mu_{x1} < 5.93$ the probability of failure exceeds the allowed value $P_{f,\text{all}} = 0.0013$. Starting from $\mu_{x1} = 0$ and moving along the line $\mu_{p1} = 3.4$ yield an increase in $P_f$ to a maximum of 0.43 % when the elliptic constraint $g(\mathbf{z})$ is approached. The probability of failure decreases thereafter further away from the constraint.

### 3.2.2 Uncorrelated normally distributed random variables with varying standard deviation

In order to illustrate the case where the standard deviation of the design variables is allowed to vary in RBDO applica-tions, assume that $\sigma_{x1}$ in (37) is not given. The results in the previous subsection is therefore modified according to

$$\mathbf{S}_{\mathbf{z},eq} = \begin{bmatrix} \sigma_{x1} & 0 \\ 0 & \sigma_{p1} \end{bmatrix} = \begin{bmatrix} \sigma_{x1} & 0 \\ 0 & 0.3 \end{bmatrix}, \quad (47)$$

and

$$\begin{cases} \mathbf{A}' = \begin{bmatrix} 0.0417\sigma_{x1}^2 & 0.0075\sigma_{x1} \\ 0.0075\sigma_{x1} & 0.0037 \end{bmatrix} \\ \mathbf{k}' = \begin{bmatrix} \sigma_{x1}(0.0833\mu_{x1} - 0.3633) \\ 0.015\mu_{x1} + 0.045 \end{bmatrix} \\ c' = 0.0417\mu_{x1}^2 - 0.3633\mu_{x1} + 1.0617 \end{cases} . \quad (48)$$

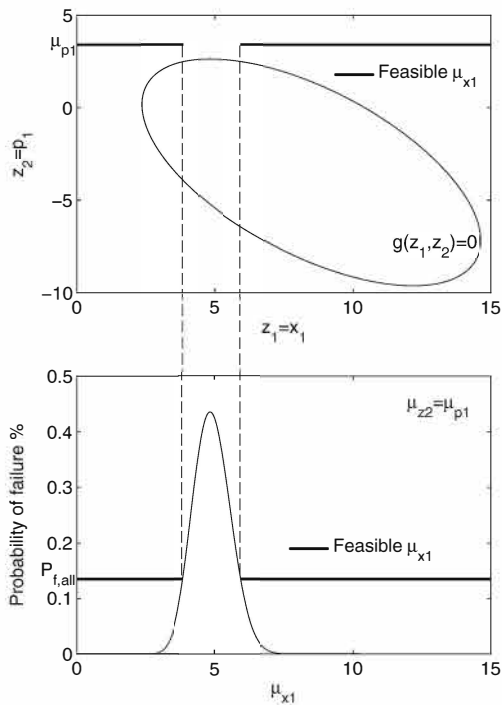

**Fig. 6** Deterministic constraint $g(x_1, p_1)$ and probability of failure as a function of mean value of standard deviation $\mu_{x1}$ for $\mu_{p1} = 3.4$

This yield the following eigenvalues

$$\gamma_i = \frac{\sigma_{x1}^2}{48} + \frac{3}{1600} \pm \frac{\sqrt{10000\sigma_{x1}^4 - 504\sigma_{x1}^2 + 81}}{4800} > 0. \quad (49)$$

Since (48) and (49) depends both on $\mu_{x1}$ and $\sigma_{x1}$, $\overline{\mathbf{k}}$ and $m_r$ in (41) and (42) will also be functions of both $\mu_{x1}$ and $\sigma_{x1}$. The probabilistic constraint in (46) therefore becomes a function of both the mean value and the standard deviation of the design variables $x$, i.e. $g^* = g^*(\mu_{x1}, \sigma_{x1})$.

### 3.2.3 Correlated arbitrarily distributed random variables with varying standard deviation

In order to demonstrate the more general case of correlated arbitrarily distributed random variables with varying standard deviation, assume that $x_1$ is log-normally distributed and that $p_1$ is normally distributed with correlation coefficient 0.5. Analogous to the previous sections, the standard deviation and mean value of design parameters are $\sigma_{p1} = 0.3$ and $\mu_{p1} = 3.4$. The matrix of correlation coefficients in (11) is given by

$$\mathbf{C} = \begin{bmatrix} 1 & 0.5 \\ 0.5 & 1 \end{bmatrix}, \quad (50)$$

and the orthogonal transformation matrix with columns consisting of the normalized eigenvectors of $\mathbf{C}$, see (12), is

$$\mathbf{T} = \frac{1}{\sqrt{2}} \begin{bmatrix} -1 & 1 \\ 1 & 1 \end{bmatrix}. \quad (51)$$

The diagonal matrix consisting of the square root of the corresponding eigenvalues of $\mathbf{C}$ according to (13), is given by

$$\mathbf{D} = \frac{1}{\sqrt{2}} \begin{bmatrix} 1 & 0 \\ 0 & \sqrt{3} \end{bmatrix}. \quad (52)$$

The probability density function (pdf) and cumulative probability density function (cdf) of the log-normally distributed variable $x$ are given by

$$\begin{cases} f_{x1}(\mu_{x1}, \sigma_{x1}) = \frac{1}{\mu_{x1}\sigma_{x1}\sqrt{2\pi}} e^{-\frac{(\ln\mu_{x1}-\mu_{x1})^2}{2\sigma_{x1}^2}} \\ F_{x1}(\mu_{x1}, \sigma_{x1}) = \frac{1}{2} + \frac{1}{2}\mathrm{erf}\left[\frac{\ln\mu_{x1}-\mu_{x1}}{\sqrt{2}\sigma_{x1}}\right] \end{cases} . \quad (53)$$

The equivalent mean value and standard deviation according to Eqn. (15) are given by

$$\begin{cases} \sigma_{x1,eq} = \frac{\phi(\Phi^{-1}[F_{x1}(\mu_{x1}, \sigma_{x1})])}{f_{x1}(\mu_{x1}, \sigma_{x1})} \\ \mu_{x1,eq} = \mu_{zi} - \Phi^{-1}[F_{x1}(\mu_{x1}, \sigma_{x1})]\sigma_{x1,eq} \end{cases} . \quad (54)$$

The vector of equivalent means $\boldsymbol{\mu}_{\mathbf{z},eq}$ and the matrix of equivalent standard deviations $\mathbf{S}_{\mathbf{z},eq}$ according to (14) are given by

$$\boldsymbol{\mu}_{\mathbf{z},eq} = \begin{bmatrix} \mu_{x1,eq} & \mu_{p1} \end{bmatrix}^T = \begin{bmatrix} \mu_{x1,eq} & 3.4 \end{bmatrix}^T \quad (55)$$

and

$$\mathbf{S}_{\mathbf{z},eq} = \begin{bmatrix} \sigma_{x1,eq} & 0 \\ 0 & \sigma_{p1} \end{bmatrix} = \begin{bmatrix} \sigma_{x1,eq} & 0 \\ 0 & 0.3 \end{bmatrix}. \quad (56)$$

Therefore $\mathbf{A}'$, $\mathbf{k}'$ and $c'$ in (20) are given by

$$\begin{cases} \mathbf{A}' = \begin{bmatrix} \frac{\sigma_{x1,eq}^2}{96} - \frac{3\sigma_{x1,eq}}{800} + \frac{3}{3200} & \frac{-\sqrt{3}(100\sigma_{x1,eq}^2-9)}{9600} \\ \frac{-\sqrt{3}(100\sigma_{x1,eq}^2-9)}{9600} & \frac{\sigma_{x1,eq}^2}{32} + \frac{9\sigma_{x1,eq}}{800} + \frac{9}{3200} \end{bmatrix} \\ \mathbf{k}' = \begin{bmatrix} \frac{1}{600}(4.5\mu_{x1,eq} + 109\sigma_{x1,eq} - 25\mu_{x1,eq}\sigma_{x1,eq} + 13.5) \\ \frac{\sqrt{3}}{1200}(9\mu_{x1,eq} - 218\sigma_{x1,eq} + 50\mu_{x1,eq}\sigma_{x1,eq} + 27) \end{bmatrix} \\ c' = \mu_{x1,eq}^2/24 - 109\mu_{x1,eq}/300 + 637/600 \end{cases}, \quad (57)$$

and the eigenvalues $\lambda_i > 0$ of $\mathbf{A}'$ are computed as

$$\lambda_i = \frac{3}{1600} + \frac{3\sigma_{x1,\text{eq}}}{800} + \frac{\sigma_{x1,\text{eq}}^2}{48}$$
$$\pm \frac{1}{4800}\sqrt{\begin{array}{l}10000\sigma_{x1,\text{eq}}^4 + 3600\sigma_{x1,\text{eq}}^3 \\ +396\sigma_{x1,\text{eq}}^2 + 324\sigma_{x1,\text{eq}} + 81\end{array}} \quad (58)$$

Observe that the equivalent standard deviation $\sigma_{x1,\text{eq}}$ and the equivalent mean value $\mu_{x1,\text{eq}}$ are functions of the standard deviation $\sigma_{x1}$ and mean value $\mu_{x1}$ according to (54). The probabilistic constraint is therefore a function of both the mean value and the standard deviation of the design variables $x$, i.e. $g^* = g^*(\mu_{x1}, \sigma_{x1})$.

*3.2.4 Deterministic and random design variables*

Assume that (31) is modified as to include a deterministic design variable $d_1$ according to

$$g(d_1, x_1, p_1) = \frac{(x_1 + p_1 - 5)^2}{30} + \frac{(x_1 - p_1 - 12)^2}{120} - d_1. \quad (59)$$

As in Section 3.2.1, the variables $x_1$ and $p_1$ are assumed normally distributed and uncorrelated with mean values $\mu_{x1}$ and $\mu_{p1} = 3.4$ and standard deviations $\sigma_{x1} = 0.3$ and $\sigma_{p1} = 0.3$. Following the same steps as in Section 3.2.1, the constraint $g$ can be written as in (10) as

$$g(\mathbf{z}) = Q(\mathbf{z}) = \mathbf{z}^T \mathbf{A} \mathbf{z} + \mathbf{k}^T \mathbf{z} + c, \quad (60)$$

where

$$\mathbf{z} = \begin{bmatrix} z_1 & z_2 & z_3 \end{bmatrix}^T = \begin{bmatrix} d_1 & x_1 & p_1 \end{bmatrix}^T, \quad (61)$$

and

$$\begin{cases} \mathbf{A} = \begin{bmatrix} 0 & 0 & 0 \\ 0 & 1/24 & 1/40 \\ 0 & 1/40 & 1/24 \end{bmatrix} \\ \mathbf{k} = \begin{bmatrix} 0 & -8/15 & -2/15 \end{bmatrix}^T \\ c = 61/30 - d_1 \end{cases} \quad (62)$$

Since the variables are uncorrelated, i.e. $\mathbf{C} = \mathbf{I}$, the orthogonal transformation matrix $\mathbf{T}$ and the diagonal matrix $\mathbf{D}$ are given by

$$\mathbf{T} = \mathbf{D} = \begin{bmatrix} 1 & 0 & 0 \\ 0 & 1 & 0 \\ 0 & 0 & 1 \end{bmatrix}. \quad (63)$$

Since the variables are normally distributed, the vector of equivalent means and the matrix of equivalent standard deviations are given by

$$\boldsymbol{\mu}_{\mathbf{z},\text{eq}} = \begin{bmatrix} d_1 & \mu_{x1} & \mu_{p1} \end{bmatrix}^T \quad (64)$$

and

$$\mathbf{S}_{\mathbf{z},\text{eq}} = \begin{bmatrix} \sigma_{d1} & 0 & 0 \\ 0 & \sigma_{x1} & 0 \\ 0 & 0 & \sigma_{p1} \end{bmatrix} = \begin{bmatrix} 0 & 0 & 0 \\ 0 & 0.3 & 0 \\ 0 & 0 & 0.3 \end{bmatrix}. \quad (65)$$

Observe that the deterministic design variable are assumed random with mean value $\mu_{d1} = d_1$ and standard deviation $\sigma_{d1} = 0$. Therefore $\mathbf{A}'$, $\mathbf{k}'$ and $c'$ in (20) are given by:

$$\begin{cases} \mathbf{A}' = \begin{bmatrix} 0 & 0 & 0 \\ 0 & 0.0037 & 0.0022 \\ 0 & 0.0022 & 0.0037 \end{bmatrix} \\ \mathbf{k}' = \begin{bmatrix} 0 \\ \mu_{x1}/40 - 0.109 \\ 3\mu_{x1}/200 + 0.045 \end{bmatrix} \\ c' = \mu_{x1}^2/24 - 0.3633\mu_{x1} + 2.0167 - d_1 \end{cases} \quad (66)$$

The orthogonal transformation matrix and the corresponding eigenvalues of $\mathbf{A}'$ are given by

$$\mathbf{P} = \begin{bmatrix} 1 & 0 & 0 \\ 0 & -1/\sqrt{2} & 1/\sqrt{2} \\ 0 & 1/\sqrt{2} & 1/\sqrt{2} \end{bmatrix} \quad (67)$$

and

$$\begin{cases} \gamma_1 = \varepsilon \\ \gamma_2 = 0.0015 \\ \gamma_3 = 0.0060 \end{cases}. \quad (68)$$

Observe that the first eigenvalue $\gamma_1 = 0$ is replaced by $\gamma_1 = \varepsilon$ where $\varepsilon$ is a sufficiently small number, see Section 3.1. In this work $\varepsilon = 10^{-7}$. Therefore following (22)

$$\bar{\mathbf{k}} = \sqrt{2}\begin{bmatrix} 0 \\ 0.077 - 0.005\mu_{x1} \\ 0.02\mu_{x1} - 0.032 \end{bmatrix} \quad (69)$$

and $m_1, m_2, m_3, m_4$ and $h$ are given by the same expressions as in (42) and (43). The expression for $q_0$ and therefore $\kappa_2$ are different than in (44) and (45) since

$$q_0 = \sum_j \frac{\bar{k}_j^2}{4\gamma_j} - c' = d_1 \quad (70)$$

and $\kappa_2$ is given by (27). The probabilistic constraint is given by (46).

## 4 RBDO applications

In this Section, three problems are solved to demonstrate the accuracy and efficiency of RSSL. In the first benchmark, a widely used mathematical example for demonstrating RBDO approaches is solved. The different steps used in the proposed RSSL approach are shown and details regarding the Design of Experiments (DOE) as well as the number of function evaluations are explained. In the second problem, the advantage of the RSSL method compared to the available RBDO methods which uses FORM for reliability assessment is demonstrated and analysed. In the last example, a challenging and widely used engineering problem for comparing different RBDO approaches is solved both using constant and varying design variance. For all examples, the proposed RSSL method is compared to the double-loop PMA approach (DL/PMA) (Tu et al. 1999), the decoupled method SORA (Du and Chen 2004) and two single-loop approaches; RDS (Shan and Wang 2008) and the single-loop method (Liang et al. 2008). Furthermore, for all methods, an SQP algorithm using MatLabs built-in function "fmincon" has been used to solve the optimization problems.

## 4.1 Mathematical Benchmark: Detailed demonstration of the RSSL method

The following mathematical problem is widely used in the literature to demonstrate and benchmark new RBDO methods, see e.g. Liang et al. (2008), Yang and Gu (2004) and Shan and Wang (2008),

$$\begin{cases} \min_{\mu_{x1},\mu_{x2}} f = \mu_{x1} + \mu_{x2} \\ \text{s.t } \Pr[g_i(\mathbf{x}) \geq 0] \geq r_{di} \\ g_1(\mathbf{x}) = \frac{x_1^2 x_2}{20} - 1 \\ g_2(\mathbf{x}) = \frac{(x_1+x_2-5)^2}{30} + \frac{(x_1-x_2-12)^2}{120} - 1 \\ g_3(\mathbf{x}) = \frac{80}{x_1^2+8x_2+5} - 1 \\ 0 \leq \mu_{xj} \leq 10, \quad j = 1, 2 \\ \sigma_{x1} = \sigma_{x2} = 0.3, \quad \beta_{di} = 3 \text{ for } i = 1, 2, 3 \end{cases} \quad (71)$$

where the desired reliabilities are given by $r_{di} = \Phi(\beta_{di})$. The problem consists of two normally distributed uncorrelated design variables. In the simulations here, the constraints are regarded as implicitly given in order to emulate industrial applications where each evaluation of a deterministic constraint typically involves one or more Finite-Element (FE) analysis. As is seen for instance in Section 4.3, typically the deterministic constraints are linked to the same FE-model or system which makes RSSL particularly efficient since the same runs are used to construct all quadratic models and therefore all probabilistic constraints. This differs from MPP methods where different runs are used to construct each probabilistic constraint since the MPP of each constraint has to be found.

The initial mean values are $\boldsymbol{\mu}_{\mathbf{x},\text{in}} = \begin{bmatrix} 5 & 5 \end{bmatrix}^T$. The deterministic problem is first solved with a total of 15 function evaluation of the constraints and objective function and solution $\boldsymbol{\mu}_{\mathbf{x},\text{det}} = \begin{bmatrix} 3.11 & 2.06 \end{bmatrix}^T$. Thereafter, quadratic response surface models of the form

$$Q(\mathbf{x}) = \mathbf{x}^T \mathbf{A} \mathbf{x} + \mathbf{k}^T \mathbf{x} + c \quad (72)$$

are fitted to each constraint around the deterministic solution $\boldsymbol{\mu}_{x,\text{det}}$. In RSSL, see Section 3.1, a Box-Behnken Design (BBD) is used with an interval for the Design of Experiment (DOE) of $\pm 1.4\beta_d \boldsymbol{\sigma}_\mathbf{x}$. This interval corresponds to $\pm 1.4\beta_d$ in the standardized normal space. However for problems with 2 random variables, BBD is not defined, see Table 1. Therefore an inscribed Central Composite Design (CCD) is used, where the star points are at $\pm 1.4\beta_d \boldsymbol{\sigma}_\mathbf{x}$. The number of function evaluations for the CCD design is 9, see Table 1. All constraints $g_1$, $g_2$ and $g_3$ are assumed to be connected to the same system, see Liang et al. (2008) and Shan and Wang (2008), and a total of 9 function evaluations are needed to construct all quadratic surrogate models $Q_1$, $Q_2$ and $Q_3$. However if the deterministic constraints are assumed to be linked to different FE-models for instance, then 9 runs would be necessary to construct each surrogate model, i.e. $3 \times 9 = 27$ runs to construct all three constraints. It should also be noted that the number of function evaluation for RDS and the Single Loop method would be roughly multiplied by 3 since the numerical gradient for each constraint need to be evaluated separately.

The RBDO problem is thereafter solved using the RSSL method with solution $\boldsymbol{\mu}_{\mathbf{x},\text{prob}} = \begin{bmatrix} 3.4454 & 3.2714 \end{bmatrix}^T$. In Fig. 7a, the deterministic constraints according to (71) as well as the deterministic solution is presented. In Fig. 7b the region for the Design of Experiment (DOE) $\pm 1.4\beta_d \boldsymbol{\sigma}_\mathbf{x}$ is marked by a rectangle and the quadratic response surface models are shown along with the exact deterministic constraints. The probabilistic constraints according to the RSSL approach are also plotted together with the exact probabilistic constraints ($10^7$ Monte Carlo simulations are used for every point on the curves). The probabilistic solution is also shown in the figure. The convergence of the objective function and the number of function evaluations are presented in Fig. 7c. As can be seen a total of 15 function evaluations of the deterministic constraints and objective function are needed to solve the deterministic problem. Thereafter 9 function evaluations are used to construct the quadratic response surface models according to (72). During the solution of the RBDO problem starting from the deterministic solution, no more function evaluation of the constraints in (71) are used since the quadratic response surface models are used to construct the RSSL probabilistic constraints. As can be seen, the RSSL constraints are accurate. In Fig. 7d the convergence of the design variables are also presented and the results are compared to other methods in Table 2. As can be seen from Table 2, the number of function evaluations in RSSL is $15 + 9 = 24$, which is much lower than double-loops and decoupled methods, but higher than FORM-based single-loop methods. However, the accuracy in RSSL is the highest. All methods yield $\beta_{\text{MC1}} \approx 2.97$ which is lower than the required reliability $\beta_d = 3$. However, RSSL is the only method that accurately achieves $\beta_{\text{MC2}} = 3.00$ while all other solution yields higher values and therefore also a higher value of the objective function.

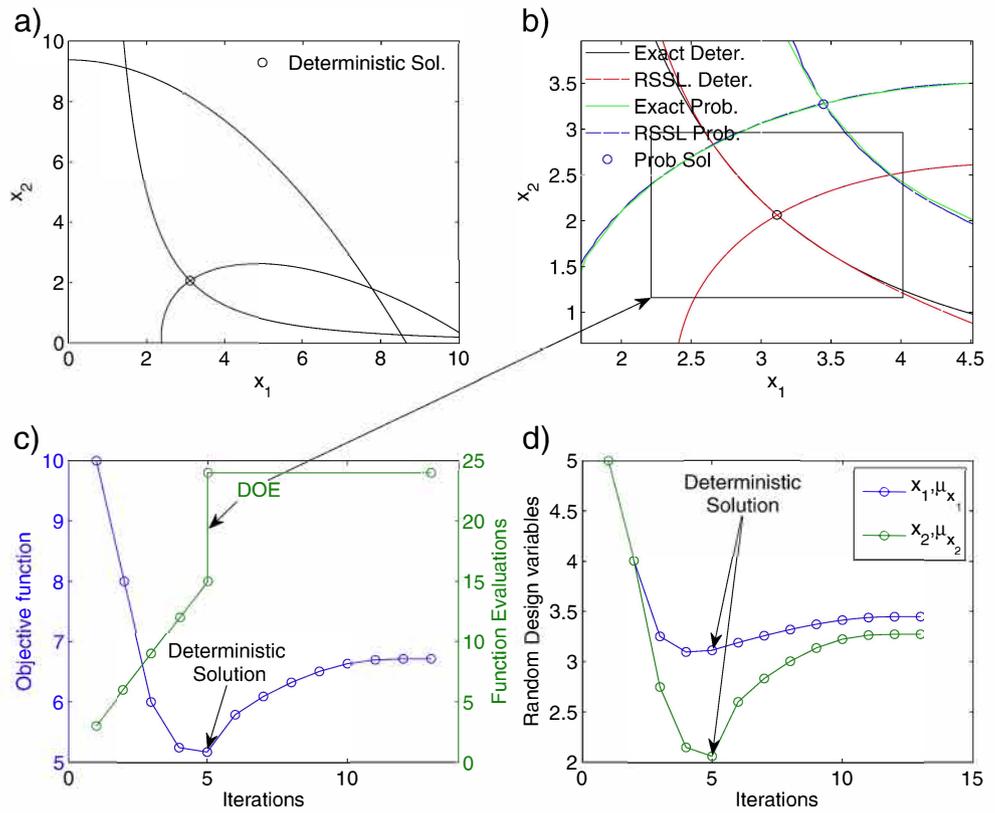

**Fig. 7** **a** Exact deterministic constraint. **b** Exact and RSSL deterministic and probabilistic constraints. **c** Convergence of objective function and number of function evaluation. **d** Convergence of design variables

The values of the generalized reliability index $\beta_{MCi}$ are computed using (2) where the probability of failure $P_{f,MCi}$ is computed using $10^7$ simulations. The accuracy in the MC computed probabilities is $10^{-4}$ with a 95 % confidence interval throughout this paper.

### 4.2 Demonstration of the limitations of FORM-based RBDO methods compared to RSSL

Consider the following mathematical problem. The limit-state functions are given by Lee et al. (2015) according to

$$\begin{cases} g_1(\mathbf{x}) = x_1^2 + 2x_1 + x_2^2 + 2x_2 - 0.5x_1x_2 - 13 \\ g_2(\mathbf{x}) = -x_1^2 - x_2^2 - x_3^2 - x_4^2 + 10x_1 + \\ \quad\quad + 12x_2 + 12x_3 + 12x_4 - 43 \end{cases}, \quad (73)$$

where $x_i \sim N(0, 1)$ for $i = 1 - 4$. A RBDO problem can be stated as

$$\begin{cases} \min_{\boldsymbol{\mu_x}} C(\boldsymbol{\mu_x}) = \mu_{x1}^2 + \mu_{x2}^2 + \mu_{x3}^2 + \mu_{x4}^2 \\ \Pr[g_i(\mathbf{x}) > 0] \leq P_{f,all}, \quad i = 1, 2 \\ -4 < \mu_{x_i} < 4 \end{cases} \quad (74)$$

The initial mean values are $\mu_{xi} = 1$ for $i = 1 - 4$. It should be noted that, although the constraints are quadratic, they are not assumed explicitly given. Therefore the same procedure as in the previous example is used here. The deterministic problem is first solved with a total of 15 function evaluations of constraints and objective function and solution $\boldsymbol{\mu}_{\mathbf{x},det} = \begin{bmatrix} 0 & 0 & 0 & 0 \end{bmatrix}^T$. Thereafter a total of 25 function

**Table 2** Test and comparison of results

|  | Deterministic | DLP/PMA | SORA | RDS | Single Loop | RSSL |
|---|---|---|---|---|---|---|
| $\mu_{x1}$ | 3.1129 | 3.4391 | 3.4385 | 3.4406 | 3.4391 | 3.4454 |
| $\mu_{x2}$ | 2.0631 | 3.2866 | 3.2871 | 3.2800 | 3.2864 | 3.2714 |
| Objective | 5.1760 | 6.7257 | 6.7256 | 6.7205 | 6.7255 | 6.7168 |
| No. of function evaluations | 15 | 1,227 | 151 | 15 | 15 | 24 |
| $\beta_{MC1}$ | **−0.0416** | **2.9744** | **2.9732** | **2.9693** | **2.96740** | **2.9725** |
| $\beta_{MC2}$ | **0.0786** | **3.0568** | **3.0591** | **3.0342** | **3.0561** | **3.0015** |

Highlighted Monte Carlo (MC) computed reliability indexes ($\beta_{MCi}$) are below the desired reliability level and therefore corresponds to violation of the probabilistic constraint

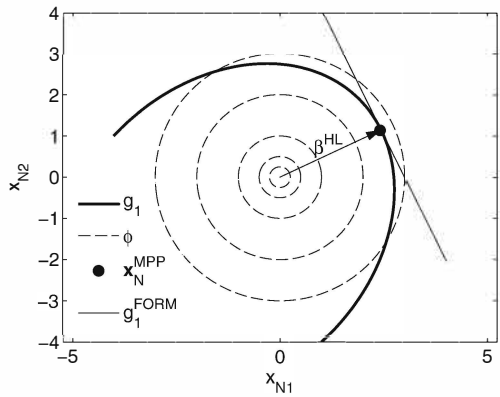

**Fig. 8** Visualisation of the nonlinearity of $g_1$ in (73) in the standard normal and uncorrelated $\mathbf{x_N}$-space at the Most Probable Point (MPP) as well as the FORM approximation of $g_1$

evaluations of the constraints in (73) are used for the BBD design. Therefore a total of $15+25 = 40$ function evaluation are used in RSSL regardless of the chosen allowed probability of failure $P_{f,all}$. As a comparison, a total of 30 function evaluation are computed for RDS and the Single-Loop Method.

This problem is simple with general quadratic models in 2 and 4 normally distributed random variables. However, the non-linearity was enough to induce high inaccuracies for RBDO methods which uses FORM for reliability assessment, see Fig. 8 and Table 3. For $P_{f,all} = 1.5$ %, the optimal solution is $\boldsymbol{\mu}_{\mathbf{x},prob} = \mathbf{0}$ and the objective function is minimized to $C = 0$. The probabilistic constraints are inactive since $P_{f1} < 1.5$ % and $P_{f2} < 1.5$ %. This is accurately given by RSSL as is seen in Table 3. However, all other methods overestimate the probability of failure with FORM, which explains the incorrect optimal solutions obtained using these methods. It leads to a higher value of the objective function and lower actual probability of failure when computed using Monte Carlo. For the more challenging case $P_{f,all} = 0.135$ % corresponding to $\beta_d = 3$, RSSL gives the solution that is closest to the optimal one. This can be seen by comparing the optimal value of the objective function which is lower in RSSL as well as comparing the values for the MC computed probabilities of failure. As can be seen in Table 3, for RSSL, the computed MC probabilities are $P_{f,MC1} = 0.1784$ % and $P_{f,MC2} = 0.1346$ % which is close to the allowed probability failure $P_{f,all} = 0.1350$ %. Although the first constraint is violated by 0.0434 percentage points, this violation is much lower than the other methods which violates the constraints by 0.2185 percentage points.

This example shows that for many non-linear problems, the accuracy of FORM-based RBDO methods, regardless of their efficiency, is questionable.

### 4.3 Reliability-Based Design for Vehicle Crashworthiness

A vehicle crashworthiness of side impact problem, see Fig. 9, has been widely used in comparing different RBDO methods, see e.g. Shan and Wang (2008), Li et al. (2010) and Du and Chen (2004). Regarding its engineering background, please refer to Gu et al. (2001). The problem is given in (75) and (76) and expressions for the deterministic constraints $g_i$ are given in Li et al. (2010), Shan and Wang (2008) and Gu et al. (2001),

$$\begin{cases} \min_{\boldsymbol{\mu}_{\mathbf{x}}} f(\boldsymbol{\mu}_{\mathbf{x}}) \\ \text{Subject to:} \\ \text{probabilistic constraints} \\ \text{Prob}\left[g_i(\mathbf{x}, \mathbf{p}) \geq 0\right] \geq r_{di} \\ \text{side constraints} \\ \boldsymbol{\mu}_{\mathbf{x}}^{Lower} \leq \boldsymbol{\mu}_{\mathbf{x}} \leq \boldsymbol{\mu}_{\mathbf{x}}^{Upper} \end{cases} \quad (75)$$

**Table 3** Test and comparison of optimal solution using different RBDO methods

|  | $P_{f,all} = 1.5$ % | | | | | $P_{f,all} = 0.1350$ % ($\beta_d = 3$) | | | | |
|---|---|---|---|---|---|---|---|---|---|---|
|  | DLP /PMA | SORA | RDS | Single Loop | RSSL | DLP /PMA | SORA | RDS | Single Loop | RSSL |
| $\mu_{x1}$ | −0.0543 | −0.0543 | −0.0539 | −0.0543 | 0 | −0.4138 | −0.4138 | −0.4138 | −0.4138 | −0.6929 |
| $\mu_{x2}$ | −0.0651 | −0.0651 | −0.0647 | −0.0651 | 0 | −0.4966 | −0.4966 | −0.4965 | −0.4966 | −0.1058 |
| $\mu_{x3}$ | −0.0651 | −0.0651 | −0.0647 | −0.0651 | 0 | −0.4966 | −0.4966 | −0.4965 | −0.4966 | −0.4331 |
| $\mu_{x4}$ | −0.0651 | −0.0651 | −0.0647 | −0.0651 | 0 | −0.4966 | −0.4966 | −0.4965 | −0.4966 | −0.4331 |
| Objective | 0.0157 | 0.0157 | 0.0155 | 0.0157 | 0 | 0.9109 | 0.9109 | 0.9106 | 0.9109 | 0.8665 |
| $P_{f,MC1}$ % | **0.5939** | **0.5939** | **0.5946** | **0.5939** | 0.7068 | **0.3535** | **0.3535** | **0.3537** | **0.3535** | **0.1784** |
| $P_{f,MC2}$ % | 1.036 | 1.036 | 1.039 | 1.036 | 1.431 | 0.08483 | 0.08483 | 0.08486 | 0.08483 | 0.1346 |

Highlighted Monte Carlo (MC) computed probability of failures ($P_{f,MCi}$) are higher than the allowed probability of failure level and therefore correspond to violation of the probabilistic constraint. The accuracy in MC computation is $10^{-4}$

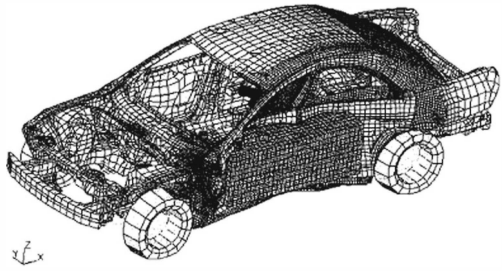

**Fig. 9** Side impact model, taken from (Gu et al. 2001)

where

$$\begin{cases} \boldsymbol{\mu}_{\mathbf{x}}^{\text{Lower}} = \begin{bmatrix} 0.5 & 0.45 & 0.5 & 0.5 & 0.875 & 0.4 & 0.4 \end{bmatrix}^T \\ \boldsymbol{\mu}_{\mathbf{x}}^{\text{Upper}} = \begin{bmatrix} 1.5 & 1.35 & 1.5 & 1.5 & 2.625 & 1.2 & 1.2 \end{bmatrix}^T \\ \text{Design parameters} \\ \boldsymbol{\mu}_{\mathbf{p}} = \begin{bmatrix} \mu_8 & \mu_9 & 0 & 0 \end{bmatrix}^T \quad \boldsymbol{\sigma}_{\mathbf{p}} = \begin{bmatrix} 0.006 & 0.006 & 10 & 10 \end{bmatrix}^T \\ \mu_8 = 0.192 \text{ or } 0.345 \quad \mu_9 = 0.192 \text{ or } 0.345 \end{cases}$$
(76)

There are 7 random design variables $x_{1-7}$ representing sizes of the structure, 2 random design parameters $p_{1-2}$ representing material properties and 2 random design parameters $p_{3-4}$ representing barrier height and barrier hitting position, respectively. All random variables are normally distributed. In this design model, $f$ is the weight of the structure which is a linear function of $\mathbf{x}$ (Gu et al. 2001). Therefore the mean of the weight, which is minimized in the problem according to (75), can be written as $f(\boldsymbol{\mu}_{\mathbf{x}})$. The deterministic constraints are (Du and Chen 2004) $g_1 = 1.0\text{kN} - F_{\text{Ab}}$ where $F_{\text{Ab}}$ is abdominal load, $g_2 = 32 \text{ mm} - D_{\text{low}}$, $g_3 = 32 \text{ mm} - D_{\text{middle}}$, $g_4 = 32 \text{ mm} - D_{\text{up}}$ where $D$'s are rib deflections (lower $D_{\text{lower}}$, middle $D_{\text{middle}}$ and upper $D_{\text{up}}$), $g_5 = 0.32 \text{ m/s} - VC_1$, $g_6 = 0.32 \text{ m/s} - VC_2$, $g_7 = 0.32 \text{ m/s} - VC_2$ where $VC$'s are viscous criteria, $g_8 = 4.01\text{kN} - F_{\text{Pubic}}$ where $F_{\text{Pubic}}$ is pubic symphysis force, $g_9 = 9.9 \text{ m/s} - V_{\text{bp}}$ where $V_{\text{bp}}$ is the velocity of B-pillar at middle point and $g_{10} = 15.69 \text{ m/s} - V_{\text{d}}$ where $V_{\text{d}}$ is the velocity of front door at B-pillar.

**Table 4** Standard deviations and initial mean values

| $\sigma_{\mathbf{x}}$ | $\mu_{\mathbf{x},\text{in}}$ |
|---|---|
| 0.03 | 1 |
| 0.03 | 0.9 |
| 0.03 | 1 |
| 0.03 | 1 |
| 0.05 | 1.75 |
| 0.03 | 0.8 |
| 0.03 | 0.8 |

All constraints are general quadratic surrogates models, presented in (Gu et al. 2001), of the form

$$g(\mathbf{z}) = Q(\mathbf{z}),$$
(77)

where $Q(\mathbf{z})$ is given by (10) and $\mathbf{z} = \begin{bmatrix} \mathbf{x}^T & \mathbf{p}^T \end{bmatrix}^T$. In this computational expensive problem, a total of only 33 CAE simulations, see Fig. 9, have been used to fit the surrogate models (Gu et al. 2001).

Using RSSL, the RBDO problem is solved in a single loop, see Fig. 3, since the quadratic surrogate models are valid on the whole region of interest. The problem is solved for two cases; constant and varying standard deviation of design variables.

*4.3.1 Constant standard deviation*

The RBDO problem is solved with constant standard deviation $\boldsymbol{\sigma}_{\mathbf{x}}$ and initial mean values $\boldsymbol{\mu}_{\mathbf{x},\text{in}}$, presented in Table 4. The convergence of the objective function is presented in Fig. 10 for the desired reliability levels $r_d = 0.9$ and $r_d = 0.9987$. The same initial mean values $\boldsymbol{\mu}_{\mathbf{x},\text{in}}$ as in Du and Chen (2004) are used resulting in an initial value of the objective function of 29.172 kg. The total number of function evaluations to convergence is compared to the ones presented by Du and Chen (2004) for SORA, double loop PMA and double loop RIA approaches, see Table 5. The number of function evaluations using the proposed RSSL approach as well as the two single loop approaches, RDS and the Single Loop method, is also presented in the table. The results presented by Du and Chen (2004) are achieved using a sequential quadratic programming (SQP) algorithm and the results presented for RSSL, RDS and the Single Loop method are achieved using MatLabs optimizer "fmincon" also with an SQP algorithm.

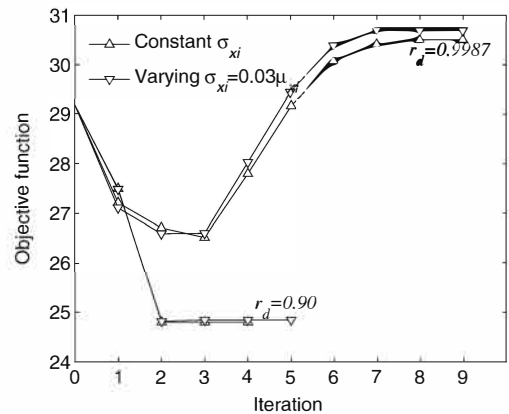

**Fig. 10** Convergence of objective function for the cases of constant and varying standard deviation for both a desired reliability level of $r_d = 0.9$ and $r_d = 0.9987$ ($\beta_d = 3$)

**Table 5** Number of function evaluation using the given quadratic surrogate models and an SQP algorithm for $r_d = 0.9$

| Method | Number of function evaluations |
| --- | --- |
| Single Loop | 32 |
| RDS | 32 |
| RSSL | 41 |
| SORA | 415 |
| Double loop PMA | 3,324 |
| Double loop RIA | 26,984 |

The results using RSSL are further compared to results achieved by other RBDO methods, see also Shan and Wang (2008), and Monte Carlo (MC) analysis with $10^9$ simulations is performed to check constraint satisfaction, see Table 6. As can be seen, RSSL is the only method that satisfies the desired reliability level of 0.9 for all constraints, whereas SORA, the Single Loop method and RDS violates the 8th probabilistic constraint by 11, 12 and 6 percentage points, respectively. For a desired reliability level of 0.99865 ($\beta_d = 3$), RSSL is still the only method satisfying the probabilistic constraints. The Single Loop Method and RDS violates the 8th and 10th constraints, and SORA violates the 8th constraints. Therefore RSSL is the only method that satisfies all constraints and minimizes the objective functions to an optimal value of 24.801 kg and 30.5009 kg for the two desired reliability levels.

*4.3.2 Varying standard deviation of design variables*

The problem is solved with a varying standard deviation given by

$$\sigma_{xi} = t_i \mu_{xi}, \quad i = 1..9 \quad (78)$$

where $t_i = 0.03$ for all standard deviations. The same initial mean values as in the constant standard deviation case are used. The computed reliability at the optimal design using Monte Carlo simulation is presented in Table 6, as can be seen all probabilistic constraints are satisfied. The convergence of the objective function is presented in Fig. 10. The total number of function evaluations for $r_d = 0.9$ is 49, which is larger but comparable to the constant standard deviation case.

**Table 6** Test and comparison of optimal solution using different RBDO methods

| | Desired reliability $r_d = 0.9$ | | | | | Desired reliability $r_d = 0.99865$ | | | | |
| --- | --- | --- | --- | --- | --- | --- | --- | --- | --- | --- |
| | SORA | Single Loop | RDS | RSSL | RSSL$_{var}$ | SORA | Single Loop | RDS | RSSL | RSSL$_{var}$ |
| $\mu_{x1}$ | 0.5000 | 0.5000 | 0.5000 | 0.5000 | 0.5000 | 0.8509 | 0.8100 | 0.8008 | 0.9336 | 0.9608 |
| $\mu_{x2}$ | 1.3091 | 1.3091 | 1.3092 | 1.3090 | 1.3143 | 1.3500 | 1.3500 | 1.3500 | 1.3500 | 1.3500 |
| $\mu_{x3}$ | 0.5000 | 0.5000 | 0.5000 | 0.5000 | 0.5000 | 0.7809 | 0.7277 | 0.7134 | 0.8994 | 0.9070 |
| $\mu_{x4}$ | 1.2942 | 1.2902 | 1.3223 | 1.3717 | 1.3737 | 1.5000 | 1.5000 | 1.5000 | 1.5000 | 1.5000 |
| $\mu_{x5}$ | 0.8750 | 0.8750 | 0.8750 | 0.8750 | 0.8750 | 0.9283 | 0.875 | 0.8750 | 0.8750 | 0.8750 |
| $\mu_{x6}$ | 1.2000 | 1.2000 | 1.2000 | 1.2000 | 1.2000 | 1.2000 | 1.2000 | 1.2000 | 1.2000 | 1.2000 |
| $\mu_{x7}$ | 0.4000 | 0.4000 | 0.4000 | 0.4000 | 0.4000 | 0.4000 | 0.4000 | 0.4000 | 0.4000 | 0.4000 |
| $\mu_{p1}$ | 0.345 | 0.345 | 0.345 | 0.345 | 0.345 | 0.345 | 0.345 | 0.345 | 0.345 | 0.345 |
| $\mu_{p2}$ | 0.192 | 0.192 | 0.192 | 0.192 | 0.192 | 0.192 | 0.192 | 0.192 | 0.192 | 0.192 |
| Objective | 24.4913 | 24.6043 | 24.6043 | 24.8010 | 24.8442 | 28.6528 | 28.6977 | 28.5526 | 30.5009 | 30.6681 |
| $r_{MC1}$ | 1.000 | 1.000 | 1.000 | 1.000 | 1.000 | 1.000 | 1.000 | 1.000 | 1.000 | 1.000 |
| $r_{MC2}$ | 0.8995 | 0.9002 | 0.9002 | 0.8999 | 0.9000 | 0.9986 | 0.9985 | 0.9986 | 0.9986 | 0.9986 |
| $r_{MC3}$ | 0.9985 | 0.9985 | 0.9985 | 0.9986 | 0.9986 | 1.000 | 1.000 | 1.000 | 1.000 | 1.000 |
| $r_{MC4}$ | 1.000 | 1.000 | 1.000 | 1.000 | 1.000 | 1.000 | 1.000 | 1.000 | 1.000 | 1.000 |
| $r_{MC5}$ | 1.000 | 1.000 | 1.000 | 1.000 | 1.000 | 1.000 | 1.000 | 1.000 | 1.000 | 1.000 |
| $r_{MC6}$ | 1.000 | 1.000 | 1.000 | 1.000 | 1.000 | 1.000 | 1.000 | 1.000 | 1.000 | 1.000 |
| $r_{MC7}$ | 1.000 | 1.000 | 1.000 | 1.000 | 1.000 | 1.000 | 1.000 | 1.000 | 1.000 | 1.000 |
| $r_{MC8}$ | **0.7879** | **0.7768** | **0.8402** | 0.9013 | 0.9000 | **0.9963** | **0.9943** | **0.9938** | 0.9985 | 0.9985 |
| $r_{MC9}$ | 0.9996 | 0.9998 | 0.9998 | 0.9999 | 0.9999 | 1.000 | 1.000 | 1.000 | 1.000 | 1.000 |
| $r_{MC10}$ | 0.9882 | 0.9928 | 0.9930 | 0.9883 | 1.000 | 0.9998 | **0.9963** | **0.9977** | 1.000 | 1.000 |

Highlighted Monte Carlo (MC) computed reliabilities ($r_{MCi}$) are below the desired reliability level and therefore corresponds to violation of the probabilistic constraint. The accuracy in MC computation is $10^{-4}$. The results using varying standard deviation are denoted RSSL$_{var}$

## 5 Discussion

### 5.1 Advantages of the RSSL-method

The proposed RSSL method is a novel approach that is not based on the concept of the Most Probable Point (MPP), as opposed to RBDO methods using FORM and SORM. The method is therefore highly efficient since locating the MPP often necessitates a large number of expensive function evaluations. This feature is of utmost importance in engineering applications where a function evaluation typically involves a large FE-run and post-processing. The cost of a function evaluation can be substantial. Furthermore, an FE-run often gives all the constraint values, see e.g. the vehicle crashworthiness problem in Section 4.3. This results in an additional strength of RSSL, since the number of function evaluations to solve the RBDO problem starting from the deterministic solution is therefore equal to the number of Design of Experiments (DoE) necessary to fit one quadratic surrogate model.

Furthermore, the available single-loop RBDO approaches, although efficient, are based on FORM which can yield substantial errors in the probability formulation. This has been shown in the second example as well as in the vehicle crashworthiness benchmark where RSSL was the only method to satisfy all probabilistic constraints. Therefore in order for an MPP method to achieve a comparable accuracy in terms of probability computation, SORM needs to be applied. However, this would involve higher computational cost and several limitations, see Section 2.2.

The proposed RSSL method can also efficiently and accurately handle problems with varying standard deviation of design variables with the same formulation, as has been shown for the vehicle crashworthiness problem.

### 5.2 Limitations of the RSSL-method

The quadratic response surface models are fitted around the deterministic solution, which can lead to errors if the distance between the MPPs of the active constraints and the deterministic solution is large or if the active constraints are highly non-linear in the region between their MPPs and the deterministic solution. For these cases, an update scheme for the surrogate model, see i.e. Wang and Wang (2014), could be used. Despite that, RSSL has been shown to be very accurate for the relatively demanding benchmarks in this paper. Nevertheless, if higher accuracy is needed, a sequential MPP-based approach with a trust region based on a zoom-pan strategy as well as analytical sensitivities of the higher-order probability of failure should be developed. This alternative approach may be more expensive than RSSL, but will further increase accuracy for critical problems since the MPP will be accurately located at each sequential iteration.

### 5.3 Accuracy and efficiency comparison of different RBDO methods

For the studied examples, it has been shown that the RSSL-method generally has higher accuracy than other RBDO methods. The efficiency of RSSL is much higher than MPP-methods and slighter lower than the single loop approaches RDS and the Single-loop Method.

It can generally be concluded that for systems that can be accurately described by quadratic response surface models, the RSSL-method yields more accurate results compared to existing RBDO methods which are based on First-Order Reliability Assessment. It should also be noted that RSSL has performed well for the mathematical problem, see Section 3.1, although the deterministic solution $\boldsymbol{\mu}_{\mathbf{x},\text{det}} = \begin{bmatrix} 3.11 & 2.06 \end{bmatrix}^T$ is not close to the probabilistic one $\boldsymbol{\mu}_{\mathbf{x},\text{prob}} = \begin{bmatrix} 3.4454 & 3.2714 \end{bmatrix}^T$. The violation of the first constraint using the DLP/PMA, SORA, RDS and Single Loop Methods is due to the use of FORM in the reliability assessment. For the proposed method, although the use of SORM in the reliability assessment yields almost exact accuracy for the second constraint, it also violates the first constraint. This is due to the residual error between the quadratic model and the exact deterministic constraint.

In general, it should also be noted that the non-linearity of the limit-state function in the standardized normal space at the Most Probable Point (MPP) is not known. It is specifically the non-linearity at the MPP which determines the accuracy or inaccuracy of FORM-based RBDO methods. Therefore, the use of the efficient RSSL method with higher-order reliability assessment is an advantage unless the limit-state function is known to be relatively linear at the MPP. This is however unlikely since the position of the MPP itself varies with the chosen mean value point, which in turn varies with the chosen desired reliability level. Therefore FORM-based RBDO methods may yield satisfactory results for a chosen reliability level with corresponding optimal mean values, but low accuracy for a different chosen value of the desired reliability level.

## 6 Conclusions

The proposed RSSL method is

- an efficient and accurate RBDO approach based on higher order probability assessment of the constraints with better accuracy than available FORM-based RBDO methods.

- an MPP-free approach with a computationally efficient formulation of the probabilistic constraints.
- a true single loop approach for problems where quadratic constraints describing the system response in the whole region of interest can be constructed.
- an RBDO method capable of solving problems both with constant and varying standard deviation of design variables.

**Acknowledgments** The financial support from Scania CV AB is gratefully acknowledged.